# SINGULARLY PERTURBED MARKOV CHAINS: LIMIT RESULTS AND APPLICATIONS

By George Yin[1] and Hanqin Zhang[2]

*Wayne State University and Academia Sinica*

This work focuses on time-inhomogeneous Markov chains with two time scales. Our motivations stem from applications in reliability and dependability, queueing networks, financial engineering and manufacturing systems, where two-time-scale scenarios naturally arise. One of the important questions is: As the rate of fluctuation of the Markov chain goes to infinity, if the limit distributions of suitably centered and scaled sequences of occupation measures exist, what can be said about the convergence rate? By combining singular perturbation techniques and probabilistic methods, this paper addresses the issue by concentrating on sequences of centered and scaled functional occupation processes. The results obtained are then applied to treat a queueing system example.

**1. Introduction.** This work focuses on time-inhomogeneous Markov chains with two time scales. Our motivations stem from applications in reliability, dependability theory, financial engineering, queueing networks and manufacturing systems, where two-time-scale scenarios naturally arise. By two time scales, we mean that the systems under consideration involve a fast varying time as well as a slowly changing one. A convenient way of formulation is to introduce a small parameter $\varepsilon > 0$. Then the fast and slow times can be represented by $t/\varepsilon$ and $t$, respectively. The characteristics of the systems associated with the slowly varying time $t$ represent the steady-state behavior, whereas those of the fast changing time $t/\varepsilon$ represent the transient behavior. One aims to reduce the complexity of the underlying systems by taking advantages of the two-time-scale formulations.

Received May 2006; revised August 2006.

[1]Supported in part by NSF Grant DMS-06-03287.

[2]Supported in part by a Distinguished Young Investigator Grant from the National Natural Sciences Foundation of China, and a grant from the Hundred Talents Program of the Chinese Academy of Sciences.

Such two-time-scale models have been used in mathematical finance by Fouque, Papanicolaou and Sircar [12], in manufacturing systems by Sethi and Zhang [28] and in control, optimization and filtering by Kushner [22]. The two-time-scale approach may also be adopted to treat reliability and dependability; see [10] and references therein for dependability models. In a wide variety of situations, one uses time-scale separation to model physical, biological and economical phenomena and to highlight the distinct rates of variations leading to a rapidly fluctuating Markov chain. To make the computation manageable and to reduce the complexity, one often uses asymptotic results to devise approximation strategies.

For example, a challenging and important problem in queueing theory is concerned with a time-varying queueing length process. Consider the queue length process of an $M_t/M_t/1/m_0$ queue, which is a single server queue with $m_0$ waiting rooms and the corresponding stochastic process is a time-inhomogeneous birth-death Markov process. The time-varying nature makes the process difficult to handle. Nevertheless, under certain circumstances, the time-dependent characteristics at time $t$ may be approximated by their quasi-stationary distribution at that time; see [16], [24] and [25]. To achieve such an approximation, one of the approaches is to model the underlying systems with the help of singular perturbation theory (see, e.g., [2]) resulting in a two-time-scale formulation leading to the use of singularly perturbed Markov chains.

Consider a continuous-time Markov chain $\alpha_\varepsilon(t)$ with a countable state space $\mathbb{N} = \{1, 2, \ldots\}$ or a finite state space $\mathcal{M} = \{1, \ldots, m_0\}$, where $\varepsilon > 0$ is a small parameter. Suppose that the infinitesimal generator of $\alpha_\varepsilon(t)$ is

$$G_\varepsilon(t) = \frac{A(t)}{\varepsilon} + B(t), \tag{1.1}$$

where $A(t)$ and $B(t)$ are themselves generators of certain Markov chains. We focus on the convergence rate of a normalized sequence of functional occupation measures of $\alpha_\varepsilon(t)$ for $t \geq 0$. When $\varepsilon$ gets smaller and smaller, the Markov chain $\alpha_\varepsilon(t)$ fluctuates more and more rapidly. Within a short period of time, the chain will reach its quasi-stationary regime. Thus, we can approximate its instantaneous behavior by its quasi-steady-state characteristics. This brings us to the problem under study with the focus on the limit behavior as $\varepsilon \to 0$. For real-world applications involving piecewise deterministic processes, we refer to [9]; for such processes having two-time-scale structures, see [25] for queueing networks, and [28] for manufacturing and production planning.

Owing to its importance in emerging applications, asymptotic properties of two-time-scale Markov chains have been studied extensively. When the Markov chain $\alpha_\varepsilon(t)$ has a finite state space, under smoothness conditions of the generators, Khasminskii, Yin and Zhang [18] constructed asymptotic expansions of probability vectors and transition probability matrices. Assuming the fast changing part of the generator consisting of several ergodic classes, Yin, Zhang and Badowski [32] showed that a suitably scaled sequence of occupation measures converges weakly to a regime-switching diffusion that is a system of diffusion processes modulated by a continuous-time Markov chain. Yin and Zhang [30] extended the results to the case in which the state space is countable. Although convergence of the scaled occupation measures has been obtained, the convergence rate has not been determined for such singularly perturbed Markov chains in [30]. Nevertheless, in many applications, it is crucial to estimate the rate of convergence of functional occupation measures to the limit process; see [20]. One of the pertinent ways of dealing with this relies on certain strong invariance principles; see [26].

There is a vast literature concerning strong invariance principles. Based on quantile transforms, Csörgö and Révész [7] obtained a strong invariance principle for partial sums of independent and identically distributed random variables. It was subsequently refined by Komlós, Major and Tusnády [19] with the best possible error bounds. Csörgö, Horváth and Steinebach [6] established strong invariance principles for renewal processes; Csörgö, Deheuvels and Horváth [5] treated random stopped sums; Kurtz [23] presented a strong approximation theorem for density-dependent Markov chains. Recently, Csáki, Csörgö, Földes and Révész [4] obtained strong approximation results for additive functions.

In this paper, our aim is to establish order of magnitude estimates of the rates of convergence of functional occupation measures to the limit processes. Our methods rely on the Skorohod representation theorem for martingales and asymptotic analysis for moments of singularly perturbed Markov chains.

The rest of the paper is organized as follows. Section 2 gives the precise formulation of the problem together with the assumptions used. Section 3 is the main part of this paper, which presents the results and the proofs. Using the Skorohod representation and asymptotic properties of two-time-scale Markov chains, we develop the rates of convergence for sequences of scaled and centered functional occupation measures. A queueing application example is discussed in Section 4. We conclude the paper in Section 5 with additional remarks.

**2. Formulation and preliminaries.** Let $(\Omega, \mathcal{F}, P)$ be a probability space. By virtue of Lemma 4.4.4 in [8], without loss of generality, we may assume that the probability space accommodates all the random variables and processes of our interests. Throughout the paper, we use $z'$ (resp. $A'$) to denote the transpose of a vector $z$ (resp. a matrix $A$), and use $\mathbb{1}$ to denote an infinite-dimensional column vector with all components being 1 [i.e., $\mathbb{1} = (1, 1, \dots)'$]. Henceforth, we write $z^i$ and $a^{ij}$ for the $i$th component of the vector $z$ and the $ij$th entry of the matrix $A = (a^{ij})$, respectively, and use a subscript to index a sequence. In addition, we often use $K$ to denote a generic positive constant with the convention $KK = K$ and $K + K = K$ used.

Let $\varepsilon > 0$ and $\alpha_\varepsilon(\cdot)$ be a time-inhomogeneous Markov chain on $(\Omega, \mathcal{F}, P)$ with countable state space $\mathbb{N}$. Suppose that the infinitesimal generator of the chain is given by (1.1). We work with a finite time horizon $t \in [0, T]$.

Recall that an infinite-dimensional matrix-valued function $G(t) = (g^{ij}(t))$ defined on $[0, T]$ is a generator of the Markov chain $\beta(t)$ if:

- $g^{ij}(\cdot)$ is Borel measurable and bounded for each $i, j \in \mathbb{N}$;
- $g^{ij}(t) \geq 0$ for all $i \neq j$; and
- $\sum_{j=1}^{\infty} g^{ij}(t) = 0$ for each $i \in \mathbb{N}$, and each bounded and Borel measurable function $f(\cdot)$ defined on $\mathbb{N}$,

$$f(\beta(t)) - \int_0^t G(s) f(\cdot)(\beta(s))\, ds \qquad \text{is a martingale,}$$

where

$$G(t) f(\cdot)(i) = \sum_{j=1}^{\infty} g^{ij}(t) f(j) \qquad \text{for each} \quad i \in \mathbb{N}.$$

Definition 2.1. A Markov chain $\beta(t)$ or its generator $G(t)$ is weakly irreducible, if the system of equations

$$\begin{aligned} \nu(t)\mathbb{1} &= 1, \\ \nu(t) G(t) &= 0, \end{aligned} \tag{2.1}$$

has a unique solution $\nu(t) = (\nu^i(t) : i \in \mathbb{N})$ with $\nu^i(t) \geq 0$ for each $i \in \mathbb{N}$, where $0 = (0, 0, \ldots)$ is an infinite-dimensional 0 vector. The unique nonnegative solution is termed a quasi-stationary distribution.

Remark 2.2. Note that an equivalent way to write (2.1) is

$$\nu(t) G_a(t) = (1 \vdots 0),$$

where $G_a(t)$ is an augmented matrix given by $G_a(t) = (\mathbb{1} \vdots G(t))$. The definition above is an extension of the weak irreducibility given in [18] in that it allows the state space of the Markov chains to be countable. Compared with the usual notion of irreducibility, it deals with time-varying generators and relaxes the usual condition by allowing some components of the quasi-stationary distribution to be equal to zero. The motivation stems from a wide variety of applications in manufacturing, production planning and queueing networks.

To proceed, let a probability row vector $p_\varepsilon(t)$ be defined by

$$\begin{aligned} p_\varepsilon(t) &= (p_\varepsilon^1(t), p_\varepsilon^2(t), \ldots, p_\varepsilon^k(t), \ldots) \\ &= (P(\alpha_\varepsilon(t)=1), P(\alpha_\varepsilon(t)=2), \ldots, P(\alpha_\varepsilon(t)=k), \ldots). \end{aligned} \tag{2.2}$$

Then $p_\varepsilon(t)$ is the solution of the forward equation

$$\begin{aligned} &\dot{p}_\varepsilon(t) = p_\varepsilon(t) G_\varepsilon(t), \\ &p_\varepsilon(0) = p(0) \qquad \text{satisfying } \; p^i(0) \geq 0 \quad \text{and} \quad \sum_{i=1}^{\infty} p^i(0) = 1. \end{aligned} \tag{2.3}$$

For notational simplicity, we have assumed the initial probability vector $p(0)$ to be independent of $\varepsilon$.

For countable-state-space Markov chains, working with infinite-dimensional vector spaces, it is natural to consider the following linear spaces:

$$\ell_1 = \left\{ v : v^i \in \mathbb{R} \text{ for each } i \in \mathbb{N} \;\text{ and }\; \sum_{i=1}^{\infty} |v^i| < \infty \right\}$$

and

$$\ell_\infty = \left\{ v : v^i \in \mathbb{R} \text{ for each } i \in \mathbb{N} \;\text{ and }\; \sup_i |v^i| < \infty \right\},$$

equipped with the norms

$$\|v\|_1 = \sum_{i=1}^{\infty} |v^i| \quad \text{and} \quad \|v\|_\infty = \sup_i |v^i|,$$

respectively (see [14], page 11). For a linear operator $A$ defined on $\ell_1$ or $\ell_\infty$, we use its induced norm $\|A\| = \sup_{\|x\|=1} \|xA\|$, where $\|\cdot\|$ denotes either the $\ell_1$ norm or the $\ell_\infty$ norm. We will use the following conditions.

(A1) For each $t \in [0, T]$, $A(t)$ is weakly irreducible. Moreover, for some $n \geq 0$, $A(\cdot) \in C^{n+2}[0, T]$ and $B(\cdot) \in C^{n+1}[0, T]$, where $C^k[0, T]$ denotes the space of functions defined on $[0, T]$ that are $k$-times continuously differentiable.

(A2) There is a $\kappa > 0$ such that for each real number $t > 0$,

$$\|\exp(A(0)t) - \mathbb{1}\nu(0)\|_\infty \leq K \exp(-\kappa t), \tag{2.4}$$

where $\nu(t) = (\nu^1(t), \nu^2(t), \ldots)$ is the quasi-stationary distribution corresponding to the generator $A(t)$.

Remark 2.3. (i) By assuming the generators to be sufficiently smooth, we can derive (uniform in the time variable $t$) asymptotic expansions of the probability vectors as well as those of the transition probability matrices.

(ii) Assumption (A2) is a Doeblin-type condition; see a similar condition in [11], page 217. A condition in a slightly different form is used in [11], page 192, in which Doob gave an illuminating discussion. Nevertheless, for our purpose of study, the current condition (2.4) appears to be more convenient. Although time-varying generators and time-inhomogeneous Markov chains are considered, (A2) is a condition only concerned with a time-invariant generator (or a constant matrix) $A(0)$. Suppose that the Markov chain generated by $A(0)$ is $\beta(t)$ and that the corresponding transition probability matrix is $P(t)=(p^{ij}(t))$. Since $P(t)=\exp(A(0)t)$, (A2) can be rephrased as a condition on the probability distribution of $\beta(t)$. Relaxing this condition using a modified norm for singularly perturbed Markov chains can be found in [1].

(iii) Note that the discussion up to now is devoted to Markov chains with a countable state space. If $\mathbb{N}$ is replaced by $\mathcal{M}$, a finite state space, the weak irreducibility of $A(t)$ in (A1) implies that (A2) holds. That is, for finite-state-space Markov chains, (A2) need not be assumed. Also, in this case, the infinite-dimensional $\nu(t)$ and $\mathbb{1}$ are replaced by finite-dimensional vectors $\nu(t)=(\nu^1(t),\dots,\nu_0^m(t))\in\mathbb{R}^{1\times m_0}$ and $\mathbb{1}=(1,\dots,1)'\in\mathbb{R}^{m_0\times 1}$, respectively.

Let

$$P_\varepsilon(t_0,t)=(p_\varepsilon^{ij}(t_0,t)),\qquad t\ge t_0, \tag{2.5}$$

be the transition matrix with

$$p_\varepsilon^{ij}(t_0,t)=P(\alpha_\varepsilon(t)=j|\alpha_\varepsilon(t_0)=i)\qquad \text{for all}\quad i,j\in\mathbb{N}.$$

Define

$$w_\varepsilon^i(s)=I_{\{\alpha_\varepsilon(s)=i\}}-\nu^i(s).$$

For $\{f(i)\}_{i=1}^\infty\in\ell_1$, define a sequence of centered functional occupation measures by

$$\begin{aligned} z_\varepsilon(t)&=\int_0^t\left(f(\alpha_\varepsilon(s))-\sum_{i=1}^\infty f(i)\nu^i(s)\right)ds\\ &=\sum_{i=1}^\infty f(i)\int_0^t w_\varepsilon^i(s)\,ds, \end{aligned} \tag{2.6}$$

and define a sequence of scaled occupation measures as

$$\xi_\varepsilon(t)=\frac{z_\varepsilon(t)}{\sqrt{\varepsilon}}. \tag{2.7}$$

Before proceeding further, we present a result on asymptotic expansions of the transition probability matrix of the process $\alpha_\varepsilon(\cdot)$ and a weak invariance principle of $\xi_\varepsilon(\cdot)$, which will be used to prove our main results. Generalizations of the above results can be found in [30], Theorems 3.7, 3.8 and 4.4.

PROPOSITION 2.4. *Assume* (A1) *and* (A2). *Then the following results hold.*

(i) *The asymptotic expansions in* (2.8) *hold uniformly in* $(t_0, t)$ *with* $0 \le t_0 \le t \le T$*:*

$$P_\varepsilon(t_0, t) = \sum_{k=0}^{n} \varepsilon^k \Phi_k(t) + \sum_{k=0}^{n} \varepsilon^k \Psi_k\biggl(t_0, \frac{t - t_0}{\varepsilon}\biggr) + O(\varepsilon^{n+1}), \tag{2.8}$$

*where* $\Phi_0(t) = \mathbb{1}\nu(t)$ *and* $\Phi_k(t)$ *and* $\Psi_k(t_0, \tau)$ $[$*with* $\tau = (t - t_0)/\varepsilon]$ *are the solutions of*

$$\Phi_k(t)A(t) = \frac{d\Phi_{k-1}(t)}{dt} - \Phi_{k-1}(t)B(t), \qquad k = 1, \dots, n,$$

$$\begin{cases} \dfrac{d\Psi_0(t_0, \tau)}{d\tau} = \Psi_0(t_0, \tau)A(t_0), \\ \Psi(t_0, t_0) = I - \Phi_0(t_0), \end{cases}$$

$$\begin{cases} \dfrac{d\Psi_k(t_0, \tau)}{d\tau} = \Psi_k(t_0, \tau)A(t_0) + R_k(t_0, \tau), \qquad k = 1, \dots, n, \\ \Psi_k(t_0, t_0) = -\Phi_k(t_0), \end{cases} \tag{2.9}$$

$$R_k(t_0, \tau) = \sum_{i=0}^{k-1} \Psi_{k-i-1}(t_0, \tau)\biggl(\frac{\tau^{i+1}}{(i+1)!}\frac{d^{i+1}A(t_0)}{dt^{i+1}} + \frac{\tau^i}{i!}\frac{d^i B(t_0)}{dt^i}\biggr),$$

*respectively, where* $I$ *is the identity matrix. Moreover, for* $k = 0, 1, \dots, n$, $\Phi_k(\cdot) \in C^{n+2-k}$, *and there exist* $K > 0$ *and* $\kappa_0 > 0$ *satisfying* $\|\Psi_k(t_0, \tau)\|_\infty \le K \exp(-\kappa_0 \tau)$.

(ii) *As* $\varepsilon \to 0$, $\xi_\varepsilon(\cdot)$ *converges weakly to a diffusion process* $\xi(\cdot)$ *such that*

$$E\xi(t) = 0 \quad \text{and} \quad E[\xi(t)]^2 = \int_0^t \sigma^2(s)\,ds \quad \text{with}$$

$$\sigma^2(s) = \sum_{i=1}^{\infty}\sum_{j=1}^{\infty} f(i)f(j)\biggl[\nu^i(s)\int_0^\infty \psi_0^{ij}(s, \tau)\,d\tau + \nu^j(s)\int_0^\infty \psi_0^{ji}(s, \tau)\,d\tau\biggr]. \tag{2.10}$$

(iii) *The following estimate holds:*

$$\sup_{0 \le t \le T}\biggl|E[\xi_\varepsilon(t)]^2 - \int_0^t \sigma^2(s)\,ds\biggr| = O(\varepsilon). \tag{2.11}$$

REMARK 2.5. Note that in deriving the weak convergence of $\xi_\varepsilon(\cdot)$ to the diffusion process, a crucial observation is that a mixing condition holds. This mixing condition can be stated as follows: For $\Delta \geq 0$, denote $\mathcal{F}_\varepsilon^{t+\Delta} = \sigma\{\alpha_\varepsilon(s) : s \geq t + \Delta\}$ and use $\mathcal{F}_{\varepsilon,t} = \{\alpha_\varepsilon(s) : s \leq t\}$. For any $\mathcal{F}_{\varepsilon,t}$-measurable $\varpi$ and $\mathcal{F}_\varepsilon^{t+\Delta}$-measurable $\varsigma$ satisfying $|\varpi| \leq 1$ and $|\varsigma| \leq 1$, there are $K > 0$ and $\kappa_1 > 0$ such that

$$
\begin{aligned}
&|E\{\varsigma|\mathcal{F}_{\varepsilon,t}\} - E\varsigma| \leq K \exp(-\kappa_1\Delta/\varepsilon) \qquad \text{a.s.}, \\
&|E[\varpi\varsigma] - E[\varpi] \cdot E[\varsigma]| \leq K \exp(-\kappa_1\Delta/\varepsilon).
\end{aligned} \tag{2.12}
$$

Condition (2.12) is a consequence of the asymptotic expansions of the transition matrices and the exponential decays of the initial layer correction terms $\Psi_k(\cdot)$ for $k = 0, \ldots, n$, owing to the Doeblin-type condition. This mixing property will be used subsequently without specific mentioning.

When the Markov chain has a finite state space $\mathcal{M}$, define the scaled occupation measures $\widetilde{\xi}_\varepsilon(t)$ similar to the countable state-space counterpart but with finite-dimensional vectors used:

$$
\widetilde{\xi}_\varepsilon(t) = \frac{1}{\sqrt{\varepsilon}}\left(\int_0^t \widetilde{w}_\varepsilon(s)\,ds\right)\widetilde{F}, \tag{2.13}
$$

where $\widetilde{w}_\varepsilon(t) = (\widetilde{w}_\varepsilon^i(t)) = (I_{\{\alpha_\varepsilon(t)=i\}} - \nu^i(t)) \in \mathbb{R}^{1\times m_0}$ and $\widetilde{F} = (f(1), \ldots, f(m_0))'$ is an arbitrary vector in $\mathbb{R}^{m_0\times 1}$. In this case, statement (i) in Proposition 2.4 remains the same, and statement (ii) is changed to: $\widetilde{\xi}_\varepsilon(\cdot)$ converges weakly to a diffusion process $\widetilde{\xi}(\cdot)$ such that

$$
E\widetilde{\xi}(t) = 0 \quad \text{and} \quad E[\widetilde{\xi}(t)]^2 = \int_0^t \widetilde{\sigma}^2(s)\,ds, \tag{2.14}
$$

where

$$
\widetilde{\sigma}^2(s) = \sum_{i=1}^{m_0}\sum_{j=1}^{m_0} f(i)f(j)\left[\nu^i(s)\int_0^\infty \psi_0^{ij}(s,r)\,dr + \nu^j(s)\int_0^\infty \psi_0^{ji}(s,r)\,dr\right];
$$

see [31], Chapters 4 and 5 for a proof.

Part (iii) indicates that the second moments of $\xi_\varepsilon(\cdot)$ and $\xi(\cdot)$ differ by a small amount. To be more precise, the error bound is of the order $O(\varepsilon)$, which is also a consequence of the asymptotic expansions.

REMARK 2.6. This paper is devoted to the convergence rate of $\xi_\varepsilon(\cdot)$ to $\xi(\cdot)$. Why is such a study needed? Although $\xi_\varepsilon(\cdot)$ [resp. $\widetilde{\xi}_\varepsilon(\cdot)$] has been shown to converge to a diffusion process, the rate of convergence is yet to be determined. Roughly, the weak convergence of the sequence of functions of scaled occupation measures indicates that $\xi_\varepsilon(\cdot)$ can be replaced by a diffusion process $\xi(\cdot)$ and the fast variations in $\xi_\varepsilon(\cdot)$ can be ignored. However, this result alone does not provide us with further information on to what extent we can ignore $\xi_\varepsilon(\cdot)$. To obtain such information and to ascertain the convergence rate are the main goals of this paper. As alluded to in the Introduction, the result will be obtained by combining singular perturbation techniques and strong invariance methods.

## 3. A sequence of functional occupation measures.

**3. A sequence of functional occupation measures.** In this section, by focusing on $\xi_\varepsilon(\cdot)$, the scaled sequence of estimation errors of the centered functional occupation measures, we develop asymptotic analysis and ascertain the rate of convergence of $\xi_\varepsilon(\cdot)$ to $\xi(\cdot)$.

THEOREM 3.1. *Assume that conditions* (A1) *and* (A2) *hold. Then there exist a constant* $\delta > 0$ *and a diffusion process* $\xi(\cdot)$ *with drift and variance satisfying* (2.10) *such that*

$$\sup_{0\le t\le T} |\xi_\varepsilon(t) - \xi(t)| =_{a.s.} o(\varepsilon^\delta). \tag{3.1}$$

REMARK 3.2. The meaning of (3.1) is

$$\lim_{\varepsilon\to 0} \frac{\xi_\varepsilon(t) - \xi(t)}{\varepsilon^\delta} = 0 \qquad \text{almost surely,}$$

and the limit holds uniformly for $t \in [0, T]$. Henceforth, this notation will be used throughout the rest of the paper.

Theorem 3.1 is stated for $\alpha_\varepsilon(t)$ having a countable state space. If the state space is $\mathcal{M}$ in lieu of $\mathbb{N}$, the statement of the theorem will be changed as follows [see Remarks 2.3(iii) and 2.5].

THEOREM 3.3. *Assume that the Markov chain has a finite state space* $\mathcal{M}$ *and* (A1) *holds. Then there exist a constant* $\delta > 0$ *and a diffusion process* $\widetilde{\xi}(\cdot)$ *with drift and variance given by* (2.14) *such that*

$$\sup_{0\le t\le T} |\widetilde{\xi}_\varepsilon(t) - \widetilde{\xi}(t)| =_{a.s.} o(\varepsilon^\delta).$$

REMARK 3.4. Theorem 3.1 should be compared with Proposition 2.4. In lieu of a weak invariance theorem, a strong invariance principle is obtained and the convergence rate is ascertained. It will be seen in the proof that we can select any positive number $\delta$ from the interval $(0, 1/4)$. In what follows, we will prove Theorem 3.1 only. The proof of Theorem 3.3 can be carried out in exactly the same way.

PROOF OF THEOREM 3.1. To facilitate the presentation, the proof is divided into three steps.

*Step* 1. Approximating $\xi_\varepsilon(t)$ by a martingale process denoted by $\widetilde{M}_\varepsilon(t)/\sqrt{\varepsilon}$ [defined in (3.5)]. The main result in this step is summarized later in Proposition 3.5.

To obtain an asymptotic martingale expression for $\xi_\varepsilon(t)$, let us define

$$\eta_\varepsilon(t) = (I_{\{\alpha_\varepsilon(t)=1\}}, I_{\{\alpha_\varepsilon(t)=2\}}, \dots, I_{\{\alpha_\varepsilon(t)=k\}}, \dots),$$
$$M_\varepsilon(t) = \eta_\varepsilon(t) - \eta_\varepsilon(0) - \int_0^t \eta_\varepsilon(s) G_\varepsilon(s)\, ds.$$

Then $M_\varepsilon(t)$ is an $\mathcal{F}_{\varepsilon,t}$-martingale with $\mathcal{F}_{\varepsilon,t} = \sigma\{\alpha_\varepsilon(s) : s \le t\}$. Using the result of [15], page 55, we can define a stochastic integral with respect to this martingale. Furthermore, the solution of the stochastic differential equation

$$d\eta_\varepsilon(t) = \eta_\varepsilon(t) G_\varepsilon(t)\, dt + dM_\varepsilon(t)$$

is given by

$$\eta_\varepsilon(t) = \eta_\varepsilon(0) P_\varepsilon(0,t) + \int_0^t dM_\varepsilon(s) P_\varepsilon(s,t),$$

where $P_\varepsilon(s,t)$ is the transition matrix given by (2.5). Noting $\Phi_0(t) = \mathbb{1}\nu(t)$ from Proposition 2.4(i), we have

$$\begin{aligned} \eta_\varepsilon(s) - \nu(s) &= \eta_\varepsilon(0)[P_\varepsilon(0,s) - \Phi_0(s)] \\ &\quad + \int_0^s dM_\varepsilon(r)\{[P_\varepsilon(r,s) - \Phi_0(s)] + \Phi_0(s)\}. \end{aligned} \tag{3.2}$$

It is easily seen that

$$\begin{aligned} f(\alpha_\varepsilon(t)) &= \sum_{i=1}^{\infty} I_{\{\alpha_\varepsilon(t)=i\}} f(i) = \eta_\varepsilon(t) F, \\ G_\varepsilon(t) f(\cdot)(\alpha_\varepsilon(t)) &= \sum_{i=1}^{\infty} I_{\{\alpha_\varepsilon(t)=i\}} G_\varepsilon(t) f(\cdot)(i) = \eta_\varepsilon(t) G_\varepsilon(t) F, \end{aligned} \tag{3.3}$$

where

$$F = (f(1), f(2), \dots, f(k), \dots)' \qquad \text{satisfying} \quad \|F\|_1 < \infty.$$

Multiplying (3.2) from the right by $F$ together with an integration over $[0,t]$ leads to

$$\int_0^t [\eta_\varepsilon(s) - \nu(s)] F\, ds - \int_0^t \eta_\varepsilon(0)[P_\varepsilon(0,s) - \Phi_0(s)] F\, ds = \widetilde{M}_\varepsilon(t), \tag{3.4}$$

where

$$\widetilde{M}_\varepsilon(t) = \int_0^t \left( \int_0^s dM_\varepsilon(r)\{[P_\varepsilon(r,s) - \Phi_0(s)] + \Phi_0(s)\} F \right) ds. \tag{3.5}$$

Note that

$$\int_0^t \int_0^s dM_\varepsilon(r)\Phi_0(s)F\,ds = \int_0^t M_\varepsilon(s)\Phi_0(s)F\,ds = 0. \tag{3.6}$$

Equation (3.6) follows from the observations

$$\sum_i I_{\{\alpha_\varepsilon(s)=i\}} = 1 \quad \text{and} \quad \sum_i I_{\{\alpha_\varepsilon(0)=i\}} = 1.$$

Thus

$$[\eta_\varepsilon(s) - \eta_\varepsilon(0)]\Phi_0(s) = [\eta_\varepsilon(s) - \eta_0(0)]\mathbb{1}\nu(s) = 0.$$

Moreover, for $u, s \geq 0$,

$$\begin{aligned} A(u)\Phi_0(s) &= A(u)\mathbb{1}\nu(s) = 0,\\ B(u)\Phi_0(s) &= B(u)\mathbb{1}\nu(s) = 0. \end{aligned}$$

As a result,

$$\begin{aligned} M_\varepsilon(s)\Phi_0(s) &= \left[\eta_\varepsilon(s) - \eta_\varepsilon(0) - \int_0^s \eta_\varepsilon(u)G_\varepsilon(u)\,du\right]\Phi_0(s)\\ &= \eta_\varepsilon(s)\Phi_0(s) - \eta_\varepsilon(0)\Phi_0(s) - \int_0^s \eta_\varepsilon(u)\left[\frac{A(u)}{\varepsilon} + B(u)\right]\Phi_0(s)\,du\\ &= 0. \end{aligned}$$

Consequently, (3.5) leads to

$$\widetilde{M}_\varepsilon(t) = \int_0^t \int_0^s dM_\varepsilon(r)[P_\varepsilon(r,s) - \Phi_0(s)]F\,ds. \tag{3.7}$$

Note that

$$\begin{aligned} \int_0^t [\eta_\varepsilon(s) - \nu(s)]F\,ds &= \int_0^t \left[f(\alpha_\varepsilon(s)) - \sum_{i=1}^{\infty} f(i)\nu^i(s)\right]ds\\ &= z_\varepsilon(t), \qquad t \in [0,T]. \end{aligned} \tag{3.8}$$

Next, define

$$X_\varepsilon(t) = \int_0^t \eta_\varepsilon(0)[P_\varepsilon(0,s) - \Phi_0(s)]F\,ds. \tag{3.9}$$

By virtue of the asymptotic expansions [see (2.8)], the almost sure boundedness of $\eta_\varepsilon(0)$ and the boundedness of $F$, we have that for any $t$,

$$\begin{aligned} \frac{|X_\varepsilon(t)|}{\sqrt{\varepsilon}} &\leq \frac{1}{\sqrt{\varepsilon}} \int_0^t \|\eta_\varepsilon(0)[P_\varepsilon(0,s) - \Phi_0(s)]\|_\infty \|F\|_1\,ds\\ &\leq \frac{K}{\sqrt{\varepsilon}} \int_0^t \exp\left(-\frac{\kappa_0 s}{\varepsilon}\right)ds\\ &= O(\sqrt{\varepsilon}) \qquad \text{a.s.} \end{aligned} \tag{3.10}$$

In view of (2.7), (3.4), (3.8) and (3.10), we have established the following result.

PROPOSITION 3.5. *Under the conditions of Theorem 3.1,*

$$\sup_{0\le t\le T}\left|\xi_\varepsilon(t)-\frac{1}{\sqrt{\varepsilon}}\widetilde{M}_\varepsilon(t)\right| =_{a.s.} O(\sqrt{\varepsilon}\,). \tag{3.11}$$

In view of the definition of $\xi_\varepsilon(t)$ and owing to (3.4) and Proposition 3.5, to prove (3.1), it suffices to show that there exists a standard Brownian motion $W(t)$ such that

$$\sup_{0\le t\le T}\left|\frac{1}{\sqrt{\varepsilon}}\widetilde{M}_\varepsilon(t)-W\left(\int_0^t\sigma^2(s)\,ds\right)\right| = o(\varepsilon^\delta) \qquad \text{a.s.} \tag{3.12}$$

for any $0<\delta<1/4$. Note that

$$\begin{aligned}
&\sup_{0\le t\le T}\left|\frac{1}{\sqrt{\varepsilon}}\widetilde{M}_\varepsilon(t)-W\left(\int_0^t\sigma^2(s)\,ds\right)\right|\\
&\quad\le \max_{1\le k\le\lfloor T/\varepsilon\rfloor}\left|\frac{1}{\sqrt{\varepsilon}}\widetilde{M}_\varepsilon(\varepsilon k)-W\left(\int_0^{\varepsilon k}\sigma^2(s)\,ds\right)\right|\\
&\qquad+\max_{1\le k\le\lfloor T/\varepsilon\rfloor}\max_{(k-1)\varepsilon\le t\le k\varepsilon}\left|\frac{1}{\sqrt{\varepsilon}}\widetilde{M}_\varepsilon(\varepsilon k)-\frac{1}{\sqrt{\varepsilon}}\widetilde{M}_\varepsilon(t)\right|\\
&\qquad+\max_{1\le k\le\lfloor T/\varepsilon\rfloor}\max_{(k-1)\varepsilon\le t\le k\varepsilon}\left|W\left(\int_0^{\varepsilon k}\sigma^2(s)\,ds\right)-W\left(\int_0^t\sigma^2(s)\,ds\right)\right|.
\end{aligned} \tag{3.13}$$

To proceed, the estimate of the terms in the second line of (3.13) is obtained in the next step, whereas the last two terms in (3.13) are dealt with in Step 3.

*Step* 2. Estimates of $|\frac{1}{\sqrt{\varepsilon}}\widetilde{M}_\varepsilon(\varepsilon k)-W(\int_0^{\varepsilon k}\sigma^2(s)\,ds)|$. The result is stated in Proposition 3.6. To prove this assertion, we need to establish Proposition 3.7 first. The proof of Proposition 3.7, in turn, will be done by proving a sequence of lemmas.

PROPOSITION 3.6. *For any* $0<\delta<1/4$,

$$\max_{1\le k\le\lfloor T/\varepsilon\rfloor}\left|\frac{1}{\sqrt{\varepsilon}}\widetilde{M}_\varepsilon(\varepsilon k)-W\left(\int_0^{\varepsilon k}\sigma^2(s)\,ds\right)\right| = o(\varepsilon^\delta) \qquad a.s. \tag{3.14}$$

PROOF. Define

$$\begin{aligned}
&\widetilde{M}_{\varepsilon,0}=\widetilde{M}_\varepsilon(0)=0,\\
&\widetilde{M}_{\varepsilon,k}=\widetilde{M}_\varepsilon(\varepsilon k) \qquad\qquad \text{for } k=1,\dots,\left\lfloor\frac{T}{\varepsilon}\right\rfloor,
\end{aligned} \tag{3.15}$$

$$Y_{\varepsilon,k} = \widetilde{M}_{\varepsilon,k} - \widetilde{M}_{\varepsilon,k-1} \qquad \text{for} \quad k = 1, \dots, \left\lfloor \frac{T}{\varepsilon} \right\rfloor.$$

From the definition of $Y_{\varepsilon,k}$, we have that

$$Y_{\varepsilon,1} = \int_0^{\varepsilon} \int_0^{s} dM_\varepsilon(r)[P_\varepsilon(r,s) - \Phi_0(s)]F\,ds \tag{3.16}$$

and

$$Y_{\varepsilon,k} = \int_{\varepsilon(k-1)}^{\varepsilon k} \int_0^{s} dM_\varepsilon(r)[P_\varepsilon(r,s) - \Phi_0(s)]F\,ds \qquad \text{for} \quad 2 \le k \le \left\lfloor \frac{T}{\varepsilon} \right\rfloor. \tag{3.17}$$

Then $\{Y_{\varepsilon,k}, 1 \le k \le \lfloor \frac{T}{\varepsilon} \rfloor\}$ is a martingale difference sequence with respect to the filtration $\{\mathcal{F}_{\varepsilon,k}, 1 \le k \le \lfloor \frac{T}{\varepsilon} \rfloor\}$, where $\mathcal{F}_{\varepsilon,k}$ denotes the $\sigma$-algebra generated by $\{Y_{\varepsilon,j} : j \le k\}$. By virtue of the martingale version of the Skorohod representation theorem (see [29], Theorem 4.3; also [27]), there exist nonnegative variables $\tau_{\varepsilon,k}$ such that

$$\begin{aligned} &\left\{ \frac{1}{\sqrt{\varepsilon}} \widetilde{M}_{\varepsilon,k} : 1 \le k \le \left\lfloor \frac{T}{\varepsilon} \right\rfloor \right\} \\ &\qquad =^{\mathbb{D}} \left\{ W(\varepsilon(\tau_{\varepsilon,1} + \cdots + \tau_{\varepsilon,k})) : 1 \le k \le \left\lfloor \frac{T}{\varepsilon} \right\rfloor \right\}, \end{aligned} \tag{3.18}$$

where $\mathbb{D}$ denotes "equal in distribution" and $W(\cdot)$ is a standard Brownian motion. Now let $\mathcal{G}_{\varepsilon,0}$ be the trivial $\sigma$-field and let $\mathcal{G}_{\varepsilon,k}$ be the $\sigma$-field generated by $\mathcal{F}_{\varepsilon,k}$ and $\{W(t) : 0 \le t \le \varepsilon \sum_{i=1}^{k} \tau_{\varepsilon,i}\}$ for $k \ge 1$. Furthermore, from Theorem 4.3 of [29], we have:

(i) each $\tau_{\varepsilon,k}$ is $\mathcal{G}_{\varepsilon,k}$-measurable;
(ii) $\varepsilon E[\tau_{\varepsilon,1}] = \varepsilon^{-1} E[Y_{\varepsilon,1}]^2$ and $\varepsilon E[\tau_{\varepsilon,k}] = \varepsilon^{-1} E[Y_{\varepsilon,k}]^2$;
(iii) for each $k$ with $2 \le k \le \lfloor \frac{T}{\varepsilon} \rfloor$,

$$\varepsilon E\{\tau_{\varepsilon,k} | \mathcal{G}_{\varepsilon,k-1}\} = \varepsilon^{-1} E\{[Y_{\varepsilon,k}]^2 | \mathcal{F}_{\varepsilon,k-1}\} \qquad \text{a.s.}$$

Using [29], Theorem 4.3 again, there is a positive constant

$$L_r = 2(8/\pi^2)^{r-1}\Gamma(r+1),$$

where $\Gamma(\cdot)$ is the familiar Gamma function such that

$$\varepsilon^r E[\tau_{\varepsilon,1}]^r \le \varepsilon^{-r} L_r E[Y_{\varepsilon,1}]^{2r}, \tag{3.19}$$

and for each $k$ with $2 \le k \le \lfloor \frac{T}{\varepsilon} \rfloor$,

$$\varepsilon^r E\{[\tau_{\varepsilon,k}]^r | \mathcal{G}_{\varepsilon,k-1}\} \le \varepsilon^{-r} L_r E\{[Y_{\varepsilon,k}]^{2r} | \mathcal{F}_{\varepsilon,k-1}\}; \tag{3.20}$$

see [13], page 269. We need to establish the following assertion.

PROPOSITION 3.7. *For any $\theta \in (0,\ 1/2)$,*

$$\max_{1\le k\le \lfloor T/\varepsilon\rfloor}\left|\varepsilon\sum_{j=1}^{k}\tau_{\varepsilon,j}-\int_0^{k\varepsilon}\sigma^2(s)\,ds\right| =_{a.s.} O(\varepsilon^{\theta}). \tag{3.21}$$

Suppose momentarily that (3.21) holds. Then

$$\begin{aligned}&\sup_{1\le k\le\lfloor T/\varepsilon\rfloor}\left|W(\varepsilon(\tau_{\varepsilon,1}+\cdots+\tau_{\varepsilon,k}))-W\left(\int_0^{k\varepsilon}\sigma^2(s)\,ds\right)\right|\\ &\qquad\le \sup_{0\le s\le\int_0^T\sigma^2(u)du}\ \sup_{0\le t\le\varepsilon^{\theta}}|W(s+t)-W(s)|.\end{aligned} \tag{3.22}$$

It follows from Theorem 1.1.1 of [8] that

$$\sup_{0\le s\le\int_0^T\sigma^2(u)du}\ \sup_{0\le t\le\varepsilon^{\theta}}|W(s+t)-W(s)| =_{a.s.} o(\sqrt{\varepsilon^{\theta}}). \tag{3.23}$$

Noting (3.18), (3.14) directly follows from (3.23). Thus it remains to prove Proposition 3.7, which is our task in the remainder of this step. □

PROOF OF PROPOSITION 3.7. Note that

$$\left\{\sum_{j=1}^{k}(\tau_{\varepsilon,j}-E\{\tau_{\varepsilon,j}|\mathcal{G}_{\varepsilon,j-1}\}),1\le k\le\left\lfloor\frac{T}{\varepsilon}\right\rfloor\right\}\qquad\text{is a martingale.}$$

Using Burkholder's inequality for martingales (see [3], Corollary 1, page 397), we obtain that for any $0<\gamma<1$,

$$\begin{aligned}&P\left(\max_{1\le k\le\lfloor T/\varepsilon\rfloor}\left|\varepsilon\sum_{j=1}^{k}(\tau_{\varepsilon,j}-E\{\tau_{\varepsilon,j}|\mathcal{G}_{\varepsilon,j-1}\})\right|>\varepsilon^{\theta}\right)\\ &\quad\le\varepsilon^{(1-\theta)(1+\gamma)}E\left(\max_{1\le k\le\lfloor T/\varepsilon\rfloor}\left|\sum_{j=1}^{k}(\tau_{\varepsilon,j}-E\{\tau_{\varepsilon,j}|\mathcal{G}_{\varepsilon,j-1}\})\right|\right)^{1+\gamma}\\ &\quad\le A_\gamma\varepsilon^{(1-\theta)(1+\gamma)}E\left(\sum_{j=1}^{\lfloor T/\varepsilon\rfloor}(\tau_{\varepsilon,j}-E\{\tau_{\varepsilon,j}|\mathcal{G}_{\varepsilon,j-1}\})^2\right)^{(1+\gamma)/2},\end{aligned} \tag{3.24}$$

where $A_\gamma$ is a constant depending only on $\gamma$. Noting the elementary inequality

$$|a+b|^r\le C_r(|a|^r+|b|^r),\qquad r>0, \tag{3.25}$$

where $C_r = 1$ or $2^{r-1}$ according to $r \le 1$ or $r \ge 1$, (3.24) and (3.25) with $r < 1$ yield that

$$\begin{aligned} & P\left(\max_{1\le k\le \lfloor T/\varepsilon\rfloor}\left|\varepsilon\sum_{j=1}^{k}(\tau_{\varepsilon,j}-E\{\tau_{\varepsilon,j}|\mathcal{G}_{\varepsilon,j-1}\})\right|>\varepsilon^{\theta}\right) \\ & \quad \le A_\gamma \varepsilon^{(1-\theta)(1+\gamma)}\sum_{j=1}^{\lfloor T/\varepsilon\rfloor} E|\tau_{\varepsilon,j}-E\{\tau_{\varepsilon,j}|\mathcal{G}_{\varepsilon,j-1}\}|^{1+\gamma} \qquad \text{for all } \gamma<1. \end{aligned} \tag{3.26}$$

Note that $\gamma < 1$ is crucial to keep the constant $A_\gamma$ independent of $\lfloor \frac{T}{\varepsilon}\rfloor$. Using Jensen's inequality for conditional expectations, (3.19) and (3.25),

$$\begin{aligned} & \varepsilon^{1+\gamma}E|\tau_{\varepsilon,j}-E\{\tau_{\varepsilon,j}|\mathcal{G}_{\varepsilon,j-1}\}|^{1+\gamma} \\ & \quad \le 2^{\gamma}\varepsilon^{1+\gamma}[E(\tau_{\varepsilon,j})^{1+\gamma}+E(E\{\tau_{\varepsilon,j}|\mathcal{G}_{\varepsilon,j-1}\})^{1+\gamma}] \\ & \quad \le 2^{1+\gamma}\varepsilon^{1+\gamma}E(\tau_{\varepsilon,j})^{1+\gamma} \\ & \quad \le 2^{1+\gamma}L_{1+\gamma}\varepsilon^{-(1+\gamma)}E[Y_{\varepsilon,j}]^{2+2\gamma}. \end{aligned} \tag{3.27}$$

We divide the rest of the proof of Proposition 3.7 into four substeps; each of the first three is presented as a lemma.

LEMMA 3.8.

$$\left(\frac{1}{\sqrt{\varepsilon}}\right)^{2+2\gamma}\sum_{k=1}^{\lfloor T/\varepsilon\rfloor}E|Y_{\varepsilon,k}|^{2+2\gamma}=O(\varepsilon^{\gamma}). \tag{3.28}$$

PROOF. Note that $\eta_\varepsilon(\cdot)$ is bounded uniformly in $t\in[0,T]$, $\varepsilon>0$, and the underlying sample point $\omega\in\Omega$ under the norm $\|\cdot\|_\infty$. That is,

$$\sup_{\omega\in\Omega}\sup_{\varepsilon>0}\sup_{0\le t\le T}\|\eta_\varepsilon(t)\|_\infty\le 1. \tag{3.29}$$

Note also that by Proposition 2.4, in particular by (2.8),

$$\|P_\varepsilon(0,s)-\Phi_0(s)\|_\infty = O\left(\varepsilon+\exp\left(-\frac{\kappa_0 s}{\varepsilon}\right)\right). \tag{3.30}$$

In view of (3.2), for any $1\le k\le \lfloor T/\varepsilon\rfloor$,

$$\begin{aligned} |Y_{\varepsilon,k}| &= \left|\int_{\varepsilon(k-1)}^{\varepsilon k}\int_0^s dM_\varepsilon(r)[P_\varepsilon(r,s)-\Phi_0(s)]F\,ds\right| \\ &= \left|\int_{\varepsilon(k-1)}^{\varepsilon k}\{[(\eta_\varepsilon(s)-\nu(s)]F-\eta_\varepsilon(0)[P_\varepsilon(0,s)-\Phi_0(s)]F\}\,ds\right| \end{aligned} \tag{3.31}$$

$$\begin{aligned}&\le \int_{\varepsilon(k-1)}^{\varepsilon k} \|\eta_\varepsilon(s)-\nu(s)\|_\infty \cdot \|F\|_1\,ds\\&\quad+\int_{\varepsilon(k-1)}^{\varepsilon k} \|\eta_\varepsilon(0)[P_\varepsilon(0,s)-\Phi_0(s)]\|_\infty \cdot \|F\|_1\,ds\\&\le K\varepsilon \qquad \text{for some } K>0.\end{aligned}$$

The constant $K$ above can be chosen to be independent of $k$ and $\omega \in \Omega$ by virtue of (3.29), (3.30) and the boundedness of $\|F\|_1$. This, in turn, implies that

$$E|Y_{\varepsilon,k}|^{2+2\gamma} \le O(\varepsilon^{2+2\gamma}).$$

Substituting the above into the left-hand side of (3.28), we arrive at

$$\begin{aligned}&\left|\left(\frac{1}{\sqrt{\varepsilon}}\right)^{2+2\gamma}\sum_{k=1}^{\lfloor T/\varepsilon\rfloor} E(Y_{\varepsilon,k})^{2+2\gamma}\right|\\&\qquad \le O(\varepsilon^{-(1+\gamma)})O(\varepsilon^{-1})O(\varepsilon^{2+2\gamma}) = O(\varepsilon^{\gamma}).\end{aligned}$$

The desired result (3.28) thus follows. □

Lemma 3.9. *For any $\theta \in (0,\ 1/2)$,*

$$\max_{1\le k\le \lfloor T/\varepsilon\rfloor}\left|\sum_{j=1}^{k}\varepsilon(\tau_{\varepsilon,j}-E\{\tau_{\varepsilon,j}|\mathcal{G}_{\varepsilon,j-1}\})\right| =_{a.s.} O(\varepsilon^{\theta}). \tag{3.32}$$

Proof. It follows from (3.24), (3.26)–(3.27) and Lemma 3.8 that

$$\begin{aligned}&P\left(\max_{1\le k\le \lfloor T/\varepsilon\rfloor}\left|\varepsilon\sum_{j=1}^{k}(\tau_{\varepsilon,j}-E\{\tau_{\varepsilon,j}|\mathcal{G}_{\varepsilon,j-1}\})\right| > \varepsilon^{\theta}\right)\\&\qquad = O(\varepsilon^{-\theta(1+\gamma)})O(\varepsilon^{\gamma}) = O(\varepsilon^{\widetilde{\theta}}) \qquad \text{with } \widetilde{\theta} = \gamma-\theta(1+\gamma).\end{aligned} \tag{3.33}$$

Note that for any $\theta \in (0,1/2)$, $\theta/(1-\theta) < 1$. Thus, for any $\theta \in (0,1/2)$, we can choose $\gamma \in (0,1)$ such that $\gamma > \theta/(1-\theta)$. This implies that $\widetilde{\theta} > 0$. Next, we pick out $\varepsilon_n = n^{-2/\widetilde{\theta}}$. Then from (3.33),

$$\sum_{n=1}^{\infty} P\left(\max_{1\le k\le \lfloor T/\varepsilon_n\rfloor}\left|\varepsilon_n\sum_{j=1}^{k}(\tau_{\varepsilon_n,j}-E\{\tau_{\varepsilon_n,j}|\mathcal{G}_{\varepsilon_n,j-1}\})\right| > \varepsilon_n^{\theta}\right) < \infty.$$

The Borel–Cantelli lemma implies that

$$\max_{1\le k\le \lfloor T/\varepsilon_n\rfloor}\left|\varepsilon_n\sum_{j=1}^{k}(\tau_{\varepsilon_n,j}-E\{\tau_{\varepsilon_n,j}|\mathcal{G}_{\varepsilon_n,j-1}\})\right| \le_{a.s.} O(\varepsilon_n^{\theta}).$$

According to the choice of $\varepsilon_n$, we have (3.32). □

LEMMA 3.10. *For any* $\theta \in (0, 1/2)$,

$$\max_{1 \le k \le \lfloor T/\varepsilon \rfloor} \left| \varepsilon \sum_{j=1}^{k} E\{\tau_{\varepsilon,j} | \mathcal{G}_{\varepsilon,j-1}\} - \int_0^{k\varepsilon} \sigma^2(s)\, ds \right| = {}_{a.s.}O(\varepsilon^\theta). \tag{3.34}$$

PROOF. Using (iii) in Step 2,

$$\varepsilon \sum_{j=1}^{k} E\{\tau_{\varepsilon,j} | \mathcal{G}_{\varepsilon,j-1}\} = \frac{1}{\varepsilon} \sum_{j=1}^{k} E\{[Y_{\varepsilon,j}]^2 | \mathcal{F}_{\varepsilon,j-1}\}. \tag{3.35}$$

Note that $\{Y_{\varepsilon,k}\}$ is an orthogonal sequence (martingale difference). Next we claim that

$$\max_{1 \le k \le \lfloor T/\varepsilon \rfloor} \left| \frac{1}{\varepsilon} \sum_{j=1}^{k} (E\{[Y_{\varepsilon,j}]^2 | \mathcal{F}_{\varepsilon,j-1}\} - E[Y_{\varepsilon,j}]^2) \right| = O(\varepsilon^\gamma). \tag{3.36}$$

This can be done by using the same techniques as in the proof of Lemma 3.8. We thus omit the details.

Similarly to $\widetilde{M}_{\varepsilon,k}$ in (3.15), define

$$\xi_{\varepsilon,0} = \xi_\varepsilon(0) = 0, \qquad \xi_{\varepsilon,k} = \xi_\varepsilon(\varepsilon k), \qquad y_{\varepsilon,k} = \xi_{\varepsilon,k} - \xi_{\varepsilon,k-1}.$$

Using the asymptotic equivalence (3.11), we have

$$\max_{1 \le k \le \lfloor T/\varepsilon \rfloor} \left| \frac{1}{\sqrt{\varepsilon}} \sum_{j=1}^{k} Y_{\varepsilon,k} - \sum_{j=1}^{k} y_{\varepsilon,k} \right| = {}_{a.s.}O(\sqrt{\varepsilon}). \tag{3.37}$$

In view of $E[\sum_{j=1}^k Y_{\varepsilon,j}]^2 = \sum_{j=1}^k E[Y_{\varepsilon,j}]^2$, (3.37) implies that

$$\max_{1 \le k \le \lfloor T/\varepsilon \rfloor} \left| \frac{1}{\varepsilon} \sum_{j=1}^{k} E[Y_{\varepsilon,j}]^2 - E\left( \sum_{j=1}^{k} y_{\varepsilon,j} \right)^2 \right| = O(\sqrt{\varepsilon}).$$

This together with (3.36) and (3.37) implies that

$$\begin{aligned} &\max_{1 \le k \le \lfloor T/\varepsilon \rfloor} \left| \varepsilon \sum_{j=1}^{k} E\{\tau_{\varepsilon,j} | \mathcal{G}_{\varepsilon,j-1}\} - \int_0^{k\varepsilon} \sigma^2(s)\, ds \right| \\ &\quad = \max_{1 \le k \le \lfloor T/\varepsilon \rfloor} \left| E\left( \sum_{j=1}^{k} y_{\varepsilon,j} \right)^2 - \int_0^{k\varepsilon} \sigma^2(s)\, ds \right| + {}_{a.s.}O(\varepsilon^\theta), \end{aligned} \tag{3.38}$$

where ${}_{a.s.}O(\varepsilon^\theta)$ indicates the order is in the sense of almost sure convergence. As a result, to derive the desired result, we need only work with $E(\sum_{j=1}^k E y_{\varepsilon,j})^2$. By telescoping and noting that $\xi_{\varepsilon,0} = 0$,

$$E\left( \sum_{j=1}^{k} y_{\varepsilon,j} \right)^2 = E\left( \sum_{j=1}^{k} (\xi_{\varepsilon,j} - \xi_{\varepsilon,j-1}) \right)^2 = E[\xi_\varepsilon(\varepsilon k)]^2. \tag{3.39}$$

Using (2.11), we arrive at

$$\max_{1\le k\le\lfloor T/\varepsilon\rfloor}\left|E[\xi_\varepsilon(\varepsilon k)]^2-\int_0^{k\varepsilon}\sigma^2(s)\,ds\right|=O(\varepsilon). \tag{3.40}$$

Combining the estimates (3.38)–(3.40) obtained thus far, the desired result (3.34) then follows. □

COMPLETION OF THE PROOF OF PROPOSITION 3.7. Combining Lemma 3.9 and Lemma 3.10 yields that for any $\theta\in(0,1/2)$,

$$\begin{aligned}\max_{1\le k\le\lfloor T/\varepsilon\rfloor}&\left|\varepsilon\sum_{\ell=1}^{k}\tau_{\varepsilon,\ell}-\int_0^{k\varepsilon}\sigma^2(s)\,ds\right|\\ &\le\max_{1\le k\le\lfloor T/\varepsilon\rfloor}\left|\sum_{j=1}^{k}\varepsilon(\tau_{\varepsilon,j}-E\{\tau_{\varepsilon,j}|\mathcal{G}_{\varepsilon,j-1}\})\right|\\ &\quad+\max_{1\le k\le\lfloor T/\varepsilon\rfloor}\left|\varepsilon\sum_{j=1}^{k}E\{\tau_{\varepsilon,j}|\mathcal{G}_{\varepsilon,j-1}\}-\int_0^{k\varepsilon}\sigma^2(s)\,ds\right|\\ &=_{a.s.}O(\varepsilon^\theta).\end{aligned}\tag{3.41}$$

Therefore, we have (3.21). □

*Step* 3. Estimates of the last two terms in (3.13). This is stated in the following proposition.

PROPOSITION 3.11. *The following estimates hold:*

$$\max_{1\le k\le\lfloor T/\varepsilon\rfloor}\max_{(k-1)\varepsilon\le t\le k\varepsilon}\left|\frac{1}{\sqrt{\varepsilon}}\widetilde{M}_\varepsilon(\varepsilon k)-\frac{1}{\sqrt{\varepsilon}}\widetilde{M}_\varepsilon(t)\right|=_{a.s.}O(\sqrt{\varepsilon}) \tag{3.42}$$

*and*

$$\begin{aligned}\max_{1\le k\le\lfloor T/\varepsilon\rfloor}\max_{(k-1)\varepsilon\le t\le k\varepsilon}&\left|W\left(\int_0^{\varepsilon k}\sigma^2(s)\right)ds-W\left(\int_0^{t}\sigma^2(s)\,ds\right)\right|\\ &=_{a.s.}o(\sqrt{\varepsilon}).\end{aligned}\tag{3.43}$$

PROOF. Similarly to (3.31),

$$\begin{aligned}\max_{1\le k\le\lfloor T/\varepsilon\rfloor}\max_{(k-1)\varepsilon\le t\le k\varepsilon}&|\widetilde{M}_\varepsilon(\varepsilon k)-\widetilde{M}_\varepsilon(t)|\\ &=\max_{1\le k\le\lfloor T/\varepsilon\rfloor}\max_{(k-1)\varepsilon\le t\le k\varepsilon}\left|\int_t^{\varepsilon k}\int_0^s dM_\varepsilon(r)[P_\varepsilon(r,s)-\Phi_0(s)]F\,ds\right|\end{aligned}$$

$$
\begin{aligned}
&= \max_{1\le k\le \lfloor T/\varepsilon\rfloor} \max_{(k-1)\varepsilon\le t\le k\varepsilon} \biggl| \int_t^{\varepsilon k} \{[\eta_\varepsilon(s)-\nu(s)]F \\
&\qquad\qquad - \eta_\varepsilon(0)[P_\varepsilon(0,s)-\Phi_0(s)]F\}\,ds \biggr| \qquad (3.44)\\
&\le \max_{1\le k\le \lfloor T/\varepsilon\rfloor} \int_{\varepsilon(k-1)}^{\varepsilon k} \|\eta_\varepsilon(s)-\nu(s)\|_\infty \cdot \|F\|_1\,ds \\
&\quad + \max_{1\le k\le \lfloor T/\varepsilon\rfloor} \int_{\varepsilon(k-1)}^{\varepsilon k} \|\eta_\varepsilon(0)[P_\varepsilon(0,s)-\Phi_0(s)]\|_\infty \cdot \|F\|_1\,ds \\
&\le K\varepsilon \qquad \text{for some} \quad K>0.
\end{aligned}
$$

Thus, (3.42) holds.

By the boundedness of $\sigma^2(s)$, it is easily seen that

$$\sup_{0\le t\le T} \left| \sum_{k=1}^{\lfloor t/\varepsilon\rfloor} \int_{\varepsilon(k-1)}^{\varepsilon k} \sigma^2(u)\,du - \int_0^t \sigma^2(u)\,du \right| \le K\varepsilon.$$

Using [8], Theorem 1.1.1 again yields (3.43). □

Completion of the proof of Theorem 3.1. It follows from (3.13), Proposition 3.6 and Proposition 3.11 that for any $\delta\in(0,1/4)$,

$$\sup_{0\le t\le T} \left| \widetilde{M}_\varepsilon(t) - W\biggl(\int_0^t \sigma^2(s)\,ds\biggr) \right| =_{a.s.} O(\varepsilon^\delta),$$

which, in view of Step 1, implies the theorem. □

Remark 3.12. In view of (3.11), (3.13) and (3.42)–(3.43), the rate of convergence of the sequence of functional occupation measures to the limit process in (3.1) is dominated by the bound of

$$\max_{1\le k\le \lfloor T/\varepsilon\rfloor} \left| \frac{1}{\sqrt{\varepsilon}} \widetilde{M}_\varepsilon(\varepsilon k) - W\biggl(\int_0^{\varepsilon k} \sigma^2(s)\,ds\biggr) \right|. \tag{3.45}$$

Since we cannot obtain that for any $\theta\ge 1/2$,

$$\max_{1\le k\le \lfloor T/\varepsilon\rfloor} \left| \varepsilon \sum_{j=1}^{k} \tau_{\varepsilon,j} - \int_0^{k\varepsilon} \sigma^2(s)\,ds \right| =_{a.s.} O(\varepsilon^\theta), \tag{3.46}$$

in view of (3.18) and Theorem 1.1.1 [8], the bound of (3.45) is $\varepsilon^\delta$ with $\delta\in(0,1/4)$.

Remark 3.13. Equation (3.21) is the key to obtain our convergence rate. This equation is based on inequality (3.20), which is given by (3.18).

The relationship (3.18) is an application of Strassen's theorem. Note that due to the time-scale separation, the small parameter $\varepsilon > 0$ is embedded in the processes, which is nonstandard. The new twist is on the use of the scaling specific to the Markov chains.

**4. Applications to queueing processes.** This section demonstrates how the results obtained can be applied to queueing problems, illustrates how the asymptotic analysis can help us to gain insight, and provides guidelines for treating a class of time-dependent queueing models. One of the main ideas is: If the rate of change is slow enough, we can approximate the time-inhomogeneous behavior by that of quasi-stationary characteristics leading to a substantial reduction of complexity.

Consider an $M_t/M_t/1/m_0$ queue with $m_0$ waiting rooms and the first-in first-out service discipline. Let $\alpha(t)$ represent the number of customers in the system at time $t$. The arrival rates are $\lambda(t)\lambda_i$ for $0 \le i \le m_0$, and the service rates are $\mu(t)\mu_i$ for $1 \le i \le m_0$, where $\lambda(t)$, $\mu(t)$, $\lambda_i$ and $\mu_i$ are all nonnegative. Then $\alpha(t)$ is a Markov process with the generator $G(t)$ given by

$$\begin{bmatrix} -\lambda(t)\lambda_0 & \lambda(t)\lambda_0 & & & \\ \mu(t)\mu_1 & -[\lambda(t)\lambda_1+\mu(t)\mu_1] & \lambda(t)\lambda_1 & & \\ & \ddots & \ddots\ \ddots & & \\ & & \mu(t)\mu_{m_0-1} & -[\lambda(t)\lambda_{m_0-1}+\mu(t)\mu_{m_0-1}] & \lambda(t)\lambda_{m_0-1} \\ & & & \mu(t)\mu_{m_0} & -\mu(t)\mu_{m_0} \end{bmatrix}.$$

For an initial time $t_0 \in [0,T]$, let $P(t_0,t)$ with $t > t_0$ be the transition matrix $(p^{ij}(t_0,t))$ with

$$p^{ij}(t_0,t) = P(\alpha(t)=j|\alpha(t_0)=i) \qquad \text{for all } 0 \le i,j \le m_0.$$

Then we have the following system of birth–death equations:

$$\frac{d}{dt}P(t_0,t) = P(t_0,t)G(t). \tag{4.1}$$

(In the above and henceforth, we add 0 to the state space to include the possibility that the queue might be empty, so $\mathcal{M} = \{0,\dots,m_0\}$. All the results obtained in the previous sections carry over.)

Assume that the rate of change of the generator $G(t)$ varies very slowly in time so that the process $P(t_0,t)$ can achieve equilibrium before there is any significant change in the rate. Following [24] and [25], we replace the arrival rate $\lambda(t)\lambda_i$ and service rate $\mu(t)\mu_i$ by $\lambda(t)\lambda_i/\varepsilon$ and $\mu(t)\mu_i/\varepsilon$, respectively, where $\varepsilon > 0$ is a smaller parameter. Also index the queue length process by $\varepsilon$, and write $\alpha(t) = \alpha_\varepsilon(t)$. Then the corresponding $P_\varepsilon(t_0,t)$ satisfies the following system of equations:

$$\frac{d}{dt}P_\varepsilon(t_0,t) = \frac{1}{\varepsilon}P_\varepsilon(t_0,t)G(t). \tag{4.2}$$

Suppose that $\lambda(\cdot) \in C^2[0,T]$ and $\mu(\cdot) \in C^2[0,T]$. Then the generator $G(t)$ satisfies (A1). By Proposition 2.4, we can obtain the probability of $\alpha_\varepsilon(t) = i$ for any time $t$, that is, Corollary 5.2 of [25]. Furthermore, we consider the occupation time of $\alpha_\varepsilon(\cdot)$ in a given state during a finite interval time $[0,t]$, which is an important performance measure for the system. To this end, let

$$\Phi_0(t) = \begin{pmatrix} \nu(t) \\ \vdots \\ \nu(t) \end{pmatrix}, \qquad \nu(t) = (\nu^0(t), \dots, \nu_0^m(t)),$$

with

$$\nu^j(t) = \left(\frac{\lambda(t)}{\mu(t)}\right)^j \prod_{k=0}^{j-1} \frac{\lambda_k}{\mu_{k+1}} \Big/ \sum_{i=0}^{m_0} \left(\frac{\lambda(t)}{\mu(t)}\right)^i \prod_{k=0}^{i-1} \frac{\lambda_k}{\mu_{k+1}}, \qquad j = 0, \dots, m_0,$$

and

$$\Psi_0(t_0,t) = (\psi_0^{ij}(t_0,t))_{(m_0+1)\times(m_0+1)} = [\mathcal{I} - \Phi_0(t_0)] \exp\left(G(0)\frac{t-t_0}{\varepsilon}\right).$$

Our previous result on weak invariance reveals that $\widetilde{\xi}_\varepsilon(\cdot)$ converges weakly to a diffusion process $\widetilde{\xi}(\cdot)$, where $\widetilde{\xi}_\varepsilon(\cdot)$ and $\widetilde{\xi}(\cdot)$ are given in Remark 2.5. Thus, loosely, $\widetilde{\xi}_\varepsilon(\cdot)$ can be replaced by the diffusion process. Nevertheless, the weak convergence alone does not provide us with any information on how close the approximation is. The needed information is provided in this paper. By virtue of Theorem 3.3, for any $i$ with $0 \le i \le m_0$, by choosing $\widetilde{F} = e_i$ the standard unit vector in $\mathbb{R}^{(m_0+1)\times 1}$, we obtain

$$\int_0^t I_{\{\alpha_\varepsilon(s)=i\}}\,ds =_{a.s} \int_0^t \nu^i(s)\,ds + \sqrt{\varepsilon} W\left(\int_0^t \widetilde{\sigma}^2(s)\,ds\right) + o(\varepsilon^{1/2+\delta}),$$
$$\delta \in (0,1/4),$$

where $W(\cdot)$ is a standard Brownian motion, and

$$\widetilde{\sigma}^2(s) = 2\nu^i(s) \int_0^\infty \psi_0^{ii}(s,r)\,dr.$$

Thus when we approximate the occupation measure

$$\int_0^t I_{\{\alpha_\varepsilon(s)=i\}}\,ds \qquad \text{by} \quad \int_0^t \nu^i(s)\,ds,$$

the approximation error is a $\sqrt{\varepsilon}$ perturbation of the standard Brownian motion. Hence, the convergence rate result will allow us to further evaluate how good the approximation is.

**5. Concluding remarks.** This work has been devoted to limit results of scaled sequences of centered functional occupation measures. Under the framework of two-time-scale formulation, using the Skorohod representation and asymptotic properties of recently developed singularly perturbed Markov chains, we have obtained convergence rate theorems, which provide almost sure invariance principles leading to the weak convergence rates of sequences of the underlying occupation measures.

The proofs of results use the Skorohod representation in an essential way; the results are most suitable for the functional form of occupation measures. Several problems are of interest from theoretical as well as practical considerations.

First, in lieu of one-dimensional functional occupation measures, we may consider vector-valued occupation processes. Somewhat different techniques are required; see, for example, Kiefer [17] and Kuelbs [21] among others. Second, the current problem setup is under the stipulation that the fast changing part of the generator corresponds to a Markov chain having a single ergodic class. A generalization of this requires the consideration that the fast varying generator consists of multiple ergodic classes. Many applications lead to such Markovian models. The corresponding weak convergence result to diffusions with regime switching has been obtained in our recent work [30]; the associated convergence rate results will be equally important. However, from one ergodic class to multiple ergodic classes is not a straightforward extension. One of the difficulties is that the scaled sequence $\xi_\varepsilon(t)$ for the multiple-ergodic-class case is no longer $\phi$-mixing. The limit process is no longer a diffusion process, but rather a switching diffusion. For each fixed $t$, the limit distribution is not purely Gaussian, but rather a Gaussian mixture. The technical details need to be thoroughly investigated, and the questions require further thought and careful consideration.

Department of Mathematics
Wayne State University
Detroit, Michigan 48202
USA
E-mail: gyin@math.wayne.edu

Institute of Applied Mathematics
Academy of Mathematics and Systems Science
Academia Sinica, Beijing 100080
China
E-mail: hanqin@amt.cn.ac